\documentclass[a4paper,12pt,oneside]{article}
\usepackage{graphicx}
\usepackage{latexsym}
\usepackage[english]{babel}
\usepackage{amsmath}
\usepackage{amstext}
\usepackage{amsthm}
\usepackage{xspace}
\usepackage{amsfonts}
\usepackage{amssymb}
\makeatletter   \input psfig.sty
\providecommand{\LyX}{L\kern-.1667em\lower.25em\hbox{Y}\kern-.125emX\@}
\makeatother \newcommand{\N}{\mathbb N} 
 
\newtheorem{theorem}{Theorem} \newtheorem{lemma}[theorem]{Lemma}
\newtheorem{proposition}[theorem]{Proposition}

\newtheorem{definition}[theorem]{Definition}

\newtheorem{remark}[theorem]{Remark}

\begin{document}

\title{{\LARGE \textbf{Information, complexity and entropy:
}}\\{\LARGE \bf a new approach to theory and measurement methods}}
\author{Vieri Benci, Claudio Bonanno, Stefano Galatolo,\\Giulia
Menconi, Federico Ponchio} \date{\it Dedicated to Giovanni Prodi who
tought us to appreciate Science and Culture in Mathematics} \maketitle

\tableofcontents

\newpage

\section{Introduction}

In this paper, we present some results on information, complexity and entropy
as defined below and we discuss their relations with the Kolmogorov-Sinai
entropy which is the most important invariant of a dynamical system. These
results have the following features and motivations:

\begin{itemize}
\item  we give a new computable definition of information and complexity which
allows to give a \textit{computable} characterization of the K-S entropy;

\item  these definitions make sense even for a \textit{single} orbit and can
be measured by suitable data compression algorithms; hence they can be used in
simulations and in the analysis of experimental data;

\item  the asymptotic behavior of these quantities allows to compute not only
the Kolmogorov-Sinai entropy but also other quantities which give a measure of
the chaotic behavior of a dynamical system even in the case of null entropy.
\end{itemize}

\subsection{Information, complexity and entropy}

Information, complexity and entropy are words which in the mathematical and
physical literature are used with different meanings and sometimes even as
synonyms (e.g. the Algorithmic Information Content of a string is called
Kolmogorov complexity; the Shannon entropy sometimes is confused with Shannon
information etc.). For this reason, an intuitive definition of these notions
as they will be used in this paper will be given soon. In our approach the
notion of \textit{information} is the basic one. Given a finite string $s$
(namely a finite sequence of symbols taken in a given alphabet), the intuitive
meaning of \textit{quantity of information} $I(s)$ contained in $s$ is the
following one:

\begin{center}
$I(s)$ \textit{is the length of the smallest binary message from which you can
reconstruct} $s$.
\end{center}

In his pioneering work, Shannon defined the quantity of information as
a statistical notion using the tools of probability theory. Thus in
Shannon framework, the quantity of information which is contained in a
string depends on its context. For example the string
$^{\prime}the^{\prime}$ contains a certain information when it is
considered as a string coming from the English language. The same
string $^{\prime}the^{\prime}$ contains much more Shannon information
when it is considered as a string coming from the Italian language
because it is much rarer in the Italian language. Roughly speaking,
the Shannon information of a string is the absolute value of the
logarithm of its probability.

However there are measures of information which depend intrinsically on the
string and not on its probability within a given context. We will adopt this
point of view. An example of these measures of information is the Algorithmic
Information Content ($AIC$). In order to define it, it is necessary to define
the partial recursive functions. We limit ourselves to give an intuitive idea
which is very close to the formal definition, so we can consider a partial
recursive function as a computer $C$ which takes a program $p$ (namely a
binary string) as an input, performs some computations and gives a string
$s=C(p)$, written in the given alphabet, as an output. The $AIC$ of a string
$s$ is defined as the shortest binary program $p$ which gives $s$ as its
output, namely
\[
I_{AIC}(s) =\min\{ |p| :C(p)=s\}
\]

where $|p|$ is the length of the string $p$. From an heuristic point of
view, the shortest program $p$ which outputs the string $s$ is a sort of
optimal encoding of $s$. The information that is necessary to reconstruct the
string is contained in the program.

Another measure of the information content of a finite string can also be
defined by a lossless data compression algorithm $Z$ which satisfies suitable
properties which will be discussed in Section \ref{fin}. We can define the
information content of the string $s$ as the length of the compressed string
$Z(s),$ namely
\[
I_{Z}\left(  s\right)  =\left|  Z(s)\right|  .
\]

The advantage of using a Compression Algorithm lies in the fact that, in this
way, the information content $I_{Z}\left(  s\right)  $ turns out to be a
computable function. For this reason we will call it Computable Information
Content $(CIC)$.

If $\omega$ is an infinite string, in general, its information is infinite;
however it is possible to define another notion: the complexity. The
complexity $K(\omega)$ of an infinite string $\omega$ is the mean
information $I$ contained in a digit of $\omega$, namely
\begin{equation}
K(\omega)=\underset{n\rightarrow\infty}{\limsup}\frac{I(\omega^{n})}{n}
\label{one}%
\end{equation}
where $\omega^{n}$ is the string obtained taking the first $n$ elements of
$\omega.$ If we equip the set of all infinite strings $\Omega$ with a
probability measure $\mu,$ then the entropy $h_{\mu}$ of $(\Omega,\mu)$ can be
defined as the expectation value of the complexity:
\begin{equation}
h_{\mu}=\int_{\Omega}K(\omega)\;d\mu\ . \label{two}%
\end{equation}

If $I\left(  \omega\right)  =I_{AIC}\left(  \omega\right)  $ or $I\left(
\omega\right)  =I_{Z}\left(  \omega\right)  ,$ under suitable assumptions on
$Z$ and $\mu,$ $h_{\mu}$ turns out to be equal to the Shannon entropy. Notice
that, in this approach, the probabilistic aspect does not appear in the
definition of information or complexity, but only in the definition of entropy.

\subsection{Dynamical systems and chaos}

Chaos, unpredictability and instability of the behavior of dynamical systems
are strongly related to the notion of information. The Kolmogorov-Sinai
entropy can be interpreted as an average measure of information that is
necessary to describe a step of the evolution of a dynamical system. The
traditional definition of the Kolmogorov-Sinai entropy is given by the methods
of probabilistic information theory. It is the translation of the Shannon
entropy into the world of dynamical systems.

We have seen that the information content of a string can be defined either
with probabilistic methods or using the $AIC$ or the $CIC$. Similarly also the
K-S entropy of a dynamical system can be defined in different ways. The
probabilistic method is the usual one, the $AIC$ method has been introduced by
Brudno \cite{brud}; the $CIC$ method has been introduced in \cite{gal4} and
\cite{ZIPPO}. So, in principle, it is possible to define the entropy of a
\emph{single} orbit of a dynamical system (which we will call, as sometimes it
has already been done in the literature, \textit{complexity of the orbit}).
There are different ways to do this (see \cite{brud}, \cite{gal2},
\cite{Gaspard}, \cite{bonanno} , \cite{gal3}). In this paper, we will
introduce a method which can be implemented in numerical simulations. Now we
will describe it briefly.

Using the usual procedure of symbolic dynamics, given a partition $\alpha$ of
the phase space of the dynamical system $\left(  X,\mu, T\right)  $, it is
possible to associate a string $\Phi_{\alpha}\left(  x\right)  $ to the orbit
having $x$ as initial condition. If $\alpha=(A_{1},\dots,A_{l}),$ then
$\Phi_{\alpha}\left(  x\right)  =(s_{0},s_{1},\dots,s_{k},\dots)$ if and only
if
\[
T^{k}x\in A_{s_{k}} \quad\forall\ k \ .
\]

If we perform an experiment, the orbit of a dynamical system can be described
only with a given degree of accuracy, described by the partition of the phase
space $X$. A more accurate measurement device corresponds to a finer partition
of $X.$ The symbolic orbit $\Phi_{\alpha}\left(  x\right)  $ is a mathematical
idealization of these measurements. We can define the complexity $K(x,\alpha)$
of the orbit with initial condition $x$ with respect to the partition $\alpha$
in the following way
\[
K(x,\alpha)=\underset{n\rightarrow\infty}\limsup\frac{I(x,\alpha,n)}{n}%
\]
where
\begin{equation}
I(x,\alpha,n):=I(\Phi_{\alpha}\left(  x\right)  ^{n}); \label{info}%
\end{equation}
here $\Phi_{\alpha}\left(  x\right)  ^{n}$ represents the first $n$ digits of
the string $\Phi_{\alpha}\left(  x\right)  .$ Letting $\alpha$ to vary among
all the \emph{computable partitions }(see Section 4.1), we set
\[
K(x)=\sup_{\alpha}K(x,\alpha)\ .
\]

The number $K(x)$ can be considered as the average amount of information
necessary to ''describe'' the orbit in the unit time when you use a
sufficiently accurate measurement device.

Notice that the complexity of each orbit $K(x)$ is defined independently of
the choice of an invariant measure. In the compact case, if $\mu$ is an
invariant measure on $X$ then $\int_{X}K(x)\;d\mu$ equals the Kolmogorov-Sinai
entropy. In other words, in an ergodic dynamical system, for almost all points
$x\in X,\;$and for suitable choice of $\alpha,\;I(x,\alpha,n)\sim h_{\mu}n.$
Notice that this result holds for a large class of Information functions $I$
as for example the $AIC$ and the $CIC.$ Thus we have obtained an alternative
way to understand of the meaning of the K-S entropy.

The above construction makes sense also for a \emph{non stationary system}.
Its average over the space $X$ is a generalization of the K-S entropy to the
non stationary case. Moreover, the asymptotic behavior of $I(x,\alpha,n)$
gives an invariant of the dynamics which is finer than the K-S entropy and is
particularly relevant when the K-S entropy is null.

It is well known that the Kolmogorov-Sinai entropy is related to the
instability of the orbits. The exact relations between the K-S entropy and the
instability of the system is given by the Ruelle-Pesin theorem (\cite{pesin}). We will recall
this theorem in the one-dimensional case. Suppose that the average rate of
separation of nearby starting orbits is exponential, namely
\[
\Delta x(n)\simeq\Delta x(0)2^{\lambda n}\;\;\;\text{for }n\ll\infty
\]
where $\Delta x(n)$ denotes the distance of these two points. The number
$\lambda$ is called Lyapunov exponent; if $\lambda>0$ the system is instable
and $\lambda$ can be considered a measure of its instability (or sensibility
with respect to the initial conditions). The Ruelle-Pesin theorem implies
that, under some regularity assumptions, $\lambda$ equals the K-S entropy.

There are chaotic dynamical systems whose entropy is null: usually they are
called weakly chaotic. Weakly chaotic dynamics arises in the study of self
organizing systems, anomalous diffusion, long range interactions and many
others. In such dynamical systems the amount of information necessary to
describe $n$ steps of an orbit is less than linear in $n$, thus the K-S
entropy is not sensitive enough to distinguish the various kinds of weakly
chaotic dynamics. Nevertheless, using the ideas we presented here, the
relation between initial data sensitivity and information content of the
orbits can be generalized to these cases.

To give an example of such a generalization, let us consider a dynamical
system $([0,1],T)$ where the transition map $T$ is
\textit{constructive\footnote{a constructive map is a map that can be defined
using a finite amount of information, see \cite{gal3}.}}, and the function
$I(x,\alpha,n)$ is defined using the $AIC$ in a slightly different way than in
Section \ref{dynsyst} (see \cite{gal3}). If the speed of separation of nearby
starting orbits goes like $\Delta x(n)\simeq\Delta x(0)f(x,n)$, then for
almost all the points $x\in\lbrack0,1]$ we have
\[
I(x,\alpha,n)\sim\log(f(x,n)).
\]

In particular, if we have power law sensitivity ( $\Delta x(n)\simeq\Delta
x(0)n^{p}$), the information content of the orbit is $I(x,\alpha,n)\sim
p\log(n)$. If we have a stretched exponential sensitivity ( $\Delta
x(n)\simeq\Delta x(0)2^{\lambda n^{p}}$, $p<1$) the information content of the
orbits will increase with the power law: $I(x,\alpha,n)\sim n^{p}.$

An example of stretched exponential is provided by the Manneville map (see
Section \ref{manne}). The Manneville map is a discrete time dynamical system
which was introduced by \cite{Manneville} as an extremely simplified model of
intermittent turbulence in fluid dynamics. The Manneville map is defined on
the unit interval to itself by $T(x)=x+x^{z}\ (\mbox{mod\ 1})$. When $z>2$ the
Manneville map is weakly chaotic and non stationary. It can be proved
\cite{Gaspard}, \cite{bonanno}, \cite{gal3} that for almost each $x$ (with
respect to the Lebesgue measure)
\begin{equation}
I_{AIC}(x,\alpha,n)\sim n^{\frac{1}{z-1}} \ . \label{ixan}%
\end{equation}

\subsection{Analysis of experimental data}

By the previous considerations, the analysis of $I(x,\alpha,n)$ gives useful
information on the underlying dynamics. Since $I(x,\alpha,n)$ can be defined
through the $CIC$, it turns out that it can be used to analyze experimental
data using a compression algorithm which satisfies the property required by
the theory and which is fast enough to analyze long strings of data. We have
implemented a particular compression algorithm we called CASToRe: Compression
Algorithm Sensitive To Regularity. CASToRe is a modification of the LZ78
algorithm. Its internal running and the heuristic motivations for such an
algorithm are described in the Appendix (see also \cite{ZIPPO}). We have used
CASToRe on the Manneville map and we have checked that the experimental
results agree with the theoretical one, namely with equation (\ref{ixan})
(Section \ref{manne}; see also \cite{ZIPPO}). Then we have used it to analyze
the behavior of $I(x,\alpha,n)$ for the logistic map at the chaos threshold
(Section \ref{logistic}, see also \cite{menconi}).

Finally, we have applied CASToRe and the CIC analysis to DNA sequences
(Section \ref{DNA}), following the ideas of \cite{Grigolini},
\cite{Grigolini96}, \cite{AcquistiBuiatti} for what concerns the study of the
randomness of symbolic strings produced by a biological source. The cited
authors performed some classical statistical techniques, so we hope that our
approach will give rise both to new answers and new questions.

\section{Information content and complexity}

\subsection{Information content of finite strings}

\label{fin} We clarify the definition of Algorithmic Information Content that
was outlined in the Introduction. For a more precise exposition see for
example \cite{Ch} and \cite{Zv}.

In the following, we will consider a finite alphabet $\mathcal{A}$,
$a=\#(\mathcal{A})$ is the cardinality of $\mathcal{A}$, and the set
$\Sigma\left(  \mathcal{A}\right)  $ of finite strings from $\mathcal{A}$,
that is $\Sigma\left(  \mathcal{A}\right)  =\bigcup_{n=1}^{\infty}%
\mathcal{A}^{n}\cup\left\{  \emptyset\right\}  $. Finally, let ${\mathcal{A}%
}^{}=\Omega_{\mathcal{A}}$ be the set of infinite strings $\psi=(\omega
_{i})_{i\in}$ with $\omega_{i} \in\mathcal{A}$ for each $i$.

Let
\[
C:\Sigma(\{0,1\})\rightarrow\Sigma(\mathcal{A})
\]
be a partial recursive function. The intuitive idea of partial recursive
function is given in the Introduction. For a formal definition we refer to any
textbook of recursion theory.

The Algorithmic Information Content $I_{AIC}(s,C)$ of a string $s$ relative to
$C$ is the length of the shortest string $p$ such that $C(p)=s$. The string
$p$ can be imagined as a program given to a computing machine and the value
$C(p)$ is the output of the computation. We require that our computing machine
is universal. Roughly speaking, a computing machine is called \emph{universal}
if it can simulate any other machine (again, for a precise definition see any
book of recursion). In particular if $U$ and $U^{\prime}$ are universal then
$I_{AIC}(s,U)\leq I_{AIC}(s,U^{\prime})+const$, where the constant depends
only on $U$ and $U^{\prime}$. This implies that, if $U$ is universal, the
complexity of $s$ with respect to $U$ depends only on $s$ up to a fixed
constant and then its asymptotic behavior does not depend on the choice of $U$.

As we said in the introduction, the shortest program which gives a string as
its output is a sort of \emph{ideal} encoding of the string. The information
which is necessary to reconstruct the string is contained in the program.

Unfortunately this coding procedure cannot be performed by any algorithm. This
is a very deep statement and, in some sense, it is equivalent to the Turing
halting problem or to the G\"odel incompleteness theorem. Then the Algorithmic
Information Content is not computable by any algorithm.

However, suppose we have some lossless (reversible) coding procedure
$Z:\Sigma(\mathcal{A})\rightarrow\Sigma(\{0,1\})$ such that from the coded
string we can reconstruct the original string (for example the data
compression algorithms that are in any personal computer). Since the coded
string contains all the information that is necessary to reconstruct the
original string, we can consider the length of the coded string as an
approximate measure of the quantity of information that is contained in the
original string.

Of course not all the coding procedures are equivalent and give the same
performances, so some care is necessary in the definition of information
content. For this reason we introduce the notion of \emph{optimality} of an
algorithm $Z$, defined by comparing its compression ratio with a statistical
measure of information: the \emph{empirical entropy}. This quantity which is
related to Shannon entropy is defined below.

Let $s$ be a finite string of length $n$. We now define $\hat{H}_{l}(s)$, the
$l^{th}$ empirical entropy of $s$. We first introduce the empirical
frequencies of a word in the string $s$: let us consider $w\in\mathcal{A}^{l}%
$, a string from the alphabet $\mathcal{A}$ with length $l$; let
$s^{(m_{1},m_{2})}\in\mathcal{A}^{m_{2}-m_{1}}$ be the string containing the
segment of $s$ starting from the $m_{1}$-th digit up to the $m_{2}$-th digit;
let
\[
\sigma(s^{(i+1,i+l)},w)=\biggl\{%
\begin{array}
[c]{cc}%
1 & if\ s^{(i+1,i+l)}=w\\
0 & otherwise
\end{array}
(\ 0 \leq i \leq n-l).
\]

The relative frequency of $w$ (the number of occurrences of the word $w$
divided by the total number of $l$-digit sub words) in $s$ is then
\[
P(s,w)=\frac1{n-l+1}\sum^{n-l}_{i=0}\sigma(s^{(i+1,i+l)},w).
\]

This can be interpreted as the ``empirical'' probability of $w$ relative to
the string $s$. Then the $l$-empirical entropy is defined by
\[
\hat{H}_{l}(s)=-\frac1{l}\sum_{w\in A^{l}}P(s,w) \log P(s,w).
\]

The quantity $l\hat{H}_{l}(s)$ is a statistical measure of the average
information content of the $l-$digit long substring of $s$.

The algorithm $Z$ is coarsely optimal if its compression ratio $|Z(s)|/|s|$
tends to be less than or equal to $\hat{H}_{k}(s)$ for each $k$.

\begin{definition}
\label{ue!} A compression algorithm $Z$ is coarsely optimal if $\forall k$
there is $f_{k}$, with $f_{k}(n)=o(n)$, such that $\forall s$ it holds
\[
\left|  {Z}(s)\right|  \leq|s|\hat{H}_{k}(s)+f_{k}(|s|)\ .
\]
\end{definition}

\begin{remark}
The universal coding algorithms LZ77 and LZ78 \cite{LZ77},\cite{LZ78}
satisfies Definition \ref{ue!}. For the proof see \cite{manzini}.
\end{remark}

However if the empirical entropy of the string is null (weak chaos)
the above definition is not satisfying (see \cite{manzini}), so we
need an algorithm having the same asymptotic behavior of the empirical
entropy. In this way even in the weakly chaotic case our algorithm
will give a meaningful measure of the information.

\begin{definition}
A compression algorithm $Z$ is optimal if there is a constant $\lambda$ such
that $\forall k$ there is a $g_{k}$ with $g_{k}(t)=o(t)$ such that $\forall s$
it holds
\[
|Z(s)|\leq\lambda|s|\hat{H}_{k}(s)+g_{k}(|Z(s)|).
\]
\end{definition}

\begin{definition}
The information content of $s$ with respect to $Z$ is defined as $I_{Z}(s)=|Z(s)|.$
\end{definition}

It is not trivial to construct an optimal algorithm. For instance the well
known Lempel-Ziv compression algorithms are not optimal (\cite{manzini}).
However the set of optimal compression algorithms is not empty. In
\cite{PerGaeta} we give an exam\-ple of optimal compression algorithm that is
similar to the the Kolmogorov frequency coding algorithm which is used also in
\cite{brud}. This compression algorithm is not of practical use because of its
computational complexity. To our knowledge the problem of finding a fast
optimal compression algorithm is still open.

\subsection{Infinite strings and complexity}

Now we show how the various definitions of information content of finite
strings can be applied to define a notion of orbit complexity for dynamical
systems. This idea has already been exploited by Brudno (\cite{brud}). However
our construction to achieve this goal is different: we use Computable
Information Content instead of the Algorithmic Information Content (as it was
done in \cite{gal4}) and computable partitions instead of open covers.

This modifications with respect to the Brudno's definition have the advantage
to give a quantity which is the limit of computable function and hence it can
be used in numerical experiments.

The relations we can prove between these notions and the entropy will be
useful as a theoretical support for the interpretation of the experimental and
numerical results. As in Brudno's approach, we will obtain that in the ergodic
case the orbit complexity has a.e. the same value as the entropy. The reader
will notice that the results which we will explain in this section are
meaningful in the positive entropy case. The null entropy cases are harder to
deal with, and they present many aspects which are not yet well understood.
There are some results based on the AIC \cite{gal3}, but there are not
theoretical results based on the CIC which, in general, are mathematically
more involved. On the other hand, our definitions based on CIC make sense and
the relative quantities can be computed numerically; in section \ref{num}, we
will present some facts based on numerical experiments.

\bigskip

First, let us consider a symbolic dynamical system (a dynamical system on a
set of infinite strings). A symbolic dynamical system is given by $(\Omega
,\mu,\sigma)$, where $\Omega=\mathcal{A}^{\mathbf{N}}$, that is $\omega
\in\Omega$ implies that $\omega$ is an infinite sequence $(\omega_{i})_{i\in\N}$
of symbols in $\mathcal{A}$, $\sigma$ is the shift map
\[
\sigma((\omega_{i})_{i\in\mathbf{N}})=(\omega_{i+1})_{i\in\mathbf{N}}%
\]
and $\mu$ is an invariant probability measure on $\Omega$. A symbolic
dynamical system can be also viewed as an information source. For the purposes
of this paper the two notions can be considered equivalent.

The complexity of an infinite string $\omega$ is the average (over all the
string $\omega$) of the quantity of information which is contained in a single
digit of $\omega$ (cfr (\ref{one})). The quantity of information of a string
can be defined in different ways; by Statistics (empirical entropy), by
Computer Science (Algorithmic Information Content) or by compression
algorithms. This will give three measures of the complexity of infinite
strings; each of them presents different kind of advantages depending on the
problem to which it is applied.

\bigskip

Let $\omega$ be an infinite string in $\Omega$. Let $\omega^{n}=(\omega
_{1}\dots\omega_{n})$ be the string containing the first $n$ digits of
$\omega$.

\begin{definition}
If $\omega\in\Omega$ then the algorithmic information complexity of $\omega$
is the average information content per digit
\[
K_{AIC}(\omega)=\mathrel{\mathop{limsup}\limits_{n\rightarrow\infty}}%
\frac{I_{AIC}(\omega^{n})}{n}\ .
\]
If $Z$ is a compression algorithm we also define the computable complexity of
$\omega$ as the asymptotic compression ratio of $Z$
\[
K_{Z}(\omega)=\mathrel{\mathop{limsup}\limits_{n\rightarrow\infty}}%
\frac{|Z(\omega^{n})|}{n}\ .
\]
\end{definition}

We also define the quantity $\hat{H}(\omega)$. If $\omega$ is an infinite
string, $\hat{H}(\omega)$ is a sort of Shannon entropy of the single string.

\begin{definition}
By the definition of empirical entropy we define:
\[
\hat{H}_{l}(\omega)=\mathrel{\mathop{limsup}\limits_{n\rightarrow\infty}%
}\hat{H}_{l}(\omega^{n})
\]
and
\[
\hat{H}(\omega)=\mathrel{\mathop{lim}\limits_{l\rightarrow\infty}}%
\hat{H}_{l}(\omega).
\]
\end{definition}

The existence of this limit is proved in \cite{LZ78}.

The following proposition is a direct consequence of ergodicity (for the proof
see again \cite{LZ78}).

\begin{proposition}
\label{6} If $(\Omega,\mu,\sigma)$ is ergodic $\hat{H}(\omega)=h_{\mu}%
(\sigma)$ (where $h_{\mu}$ is the Kol\-mo\-go\-rov-Sinai entropy) for $\mu
$-almost each $\omega$.
\end{proposition}

Moreover from the definition of coarse optimality it directly follows that:

\begin{remark}
If $Z$ is coarsely optimal then for each $k$
\[
K_{Z}(\omega)\leq\hat{H}_{k}(\omega)\ ,
\]
so that
\[
K_{Z}(\omega)\leq\hat{H}(\omega)\ .
\]
\end{remark}

\begin{remark}
As it is intuitive, the compression ratio of $Z$ cannot be less than the
average information per digit of the algorithmic information (see \cite{gal4}):$$\forall
s\ K_{Z}(s)\geq K_{AIC}(s)\ .$$ 
\end{remark}

Then we have the following result.

\begin{theorem}
\label{11} If $(\Omega,\mu,T)$ is a symbolic dynamical system and $\mu$ is
ergodic, then for $\mu$-almost each $\omega$
\[
K_{Z}(\omega)=\hat{H}(\omega)=K_{AIC}(\omega)=h_{\mu}(T)\ .
\]
\end{theorem}

\emph{Proof.} 1) We have that $K_{Z}(\omega)\geq K_{AIC}(\omega)$ and then by
the Brudno theorem (\cite{brud}) $K_{Z}(\omega)\geq h_{\mu}(T)$ for $\mu
$-almost each $\omega$.

2) On the other hand, $K_{Z}(\omega)\leq\hat{H}(\omega)$ and by Proposition
\ref{6} $\hat{H} (\omega)=h_{\mu}$ for $\mu$ almost each $\omega\in\Omega$. \qed

\section{Dynamical Systems}

\label{dynsyst}

\subsection{Information and the Kolmogorov-Sinai entropy}

Now we consider a dynamical system $(X,\mu,T)$, where $X$ is a compact metric
space, $T$ is a continuous map $T:X\rightarrow X$ and $\mu$ is a Borel
probability measure on $X$ invariant for $T$. If $\alpha=\{A_{1},\dots
,A_{n}\}$ is a measurable partition of $X$ (a partition of $X$ where the sets
are measurable) then we can associate to $(X,\mu,T)$ a symbolic dynamical
system $(\Omega_{\alpha},\mu_{\alpha},\sigma)$ ( called a symbolic model of
$(X,T)$). The set $\Omega_{\alpha}$ is a subset of $\{1,\dots,n\}^{\mathbf{N}%
}$ (the space of infinite strings made of symbols from the alphabet
$\{1,\dots,n\}$). To a point $x\in X$ it is associated a string $\omega
=(\omega_{i})_{i\in\N}=\Phi_{\alpha}(x)$ defined as
\[
\Phi_{\alpha}(x)=\omega\iff\forall j\in,T^{j}(x)\in A_{\omega_{j}}%
\]
and $\Omega_{\alpha}=\mathrel{\mathop{\cup}\limits_{x\in X}}\Phi_{\alpha}(x)$.
The measure $\mu$ on $X$ induces a measure $\mu_{\alpha}$ on the associated
symbolic dynamical system. The measure is first defined on the
cylinders\footnote{We recall that $\omega^{(k,n)}=(\omega_{i})_{k\leq i\leq
n}=(\omega_{k},\omega_{k+1},\dots,\omega_{n})$.}
\[
C(\omega^{(k,n)})=\{\overline{\omega}\in\Omega_{\alpha}:\overline{\omega}%
_{i}=\omega_{i}\ for\ k\leq i\leq n-1\}
\]
by
\[
\mu_{\alpha}(C(\omega^{(k,n)}))=\mu(\cap_{k}^{n-1}T^{-i}(A_{\omega_{i}}))
\]
and then extended by the classical Kolmogorov theorem to a measure
$\mu_{\alpha}$ on $\Omega_{\alpha}$.

\begin{definition}
We define the complexity of the orbit of a point $x\in X$, with respect to the
partition $\alpha$, as
\[
K_{AIC}(x,\alpha)=K_{AIC}(\omega),
\]%
\[
K_{Z}(x,\alpha)=K_{Z}(\omega),
\]
where $\omega=\Phi_{\alpha}(x)$.
\end{definition}

\begin{theorem}
\label{14} If $Z$ is coarsely optimal, $(X,\mu,T)$ is an ergodic dynamical
system and $\alpha$ is a measurable partition of $X$, then for $\mu$-almost
all $x$
\[
K_{Z}(x,\alpha)=h_{\mu}(T,\alpha)
\]
where $h_{\mu}(T,\alpha)$ is the metric entropy of $(X,\mu,T)$ with respect to
the measurable partition $\alpha$.
\end{theorem}

\emph{Proof}. To a dynamical system $(X,\mu,T)$ and a measurable partition
$\alpha$ it is associated a symbolic dynamical system $(\Omega_{\alpha}%
,\mu_{\alpha},\sigma)$ as seen before. If $(X,\mu,T)$ is ergodic then
$(\Omega_{\alpha},\mu_{\alpha},\sigma)$ is ergodic and $h_{\mu}(T|\alpha)$ on
$X$ equals $h_{\mu_{\alpha}}(\sigma)$ on $\Omega_{\alpha}$ (see
e.g.\cite{brud}). Now by Theorem \ref{11} for almost all points in
$\Omega_{\alpha}$ $K_{Z}(\omega)=h_{\mu_{\alpha}}(\Omega_{\alpha},\sigma)$. If
we consider $Q_{\Omega_{\alpha}}:=\{\omega\in\Omega_{\alpha}:K_{Z}%
(\omega)=h_{\mu_{\alpha}}(\Omega_{\alpha},\sigma)\}$ and $Q:=\Phi_{\alpha
}^{-1}(Q_{\Omega_{\alpha}})$ we have
\[
\forall x\in Q\quad K_{Z}(x,\alpha)=K_{Z}(\Phi_{\alpha}(x))=h_{\mu_{\alpha}%
}(\sigma)=h_{\mu}(T|\alpha)\ .
\]
According to the way in which the measure $\mu_{\alpha}$ is constructed we
have $\mu(Q)=\mu_{\alpha}(Q_{\Omega_{\alpha}})=1$. \qed
\vskip0.5cm Let $\beta_{i}$ be a family of measurable partitions such that
$\mathrel{\mathop{\lim}\limits_{ i\rightarrow\infty}}\,diam(\beta_{i})=0$. If
we consider $\mathrel{\mathop{\lim\sup
}\limits_{i\rightarrow\infty}}K_{Z}(x,\beta_{i})$ we have the following

\begin{lemma}
\label{15} If $(X,\mu,T)$ is compact and ergodic, then for $\mu$-almost all
points $x\in X$, $\mathrel{\mathop{\lim\sup
}\limits_{i\rightarrow\infty}}K_{Z}(x,\beta_{i})=h_{\mu}(T)$.
\end{lemma}

\emph{Proof.} The points for which $K_{Z}(x,\beta_{i})\neq h_{\mu}(T|\beta
_{i})$ are a set of null measure for each $i$ (Theorem \ref{14}). When
excluding all these points, we exclude (for each $i$) a zero-measure set. For
all the other points we have $K_{Z}(x,\beta_{i})=h_{\mu}(T|\beta_{i})$ and
then $\mathrel{\mathop{\lim\sup
}\limits_{i\rightarrow\infty}}K_{Z}(x,\beta_{i})=\mathrel{\mathop{\lim\sup
}\limits_{i\rightarrow\infty
}}h_{\mu}(T,\beta_{i})$. Since the diameter of the partitions $\beta_{i}$
tends to 0, we have that $\mathrel{\mathop{\lim\sup
}\limits_{i\rightarrow\infty}}h_{\mu}(T,\beta_{i})=h_{\mu}(T)$ (see e.g.
\cite{Katok} page 170), and the statement is proved. 

\qed
\vskip0.5cm \noindent\textbf{Remarks.} Theorem \ref{14} and the above lemma
show that if a system has an invariant measure, its entropy can be found by
averaging the complexity of its orbits over the invariant measure. Then, as we
saw in the introduction, the entropy may be alternatively defined as the
average orbit complexity. However if we fix a single point, its orbit
complexity is not yet well defined because it depends on the choice of a
partition. It is not possible to get rid of this dependence by taking the
supremum over all partitions (as in the construction of Kolmogorov entropy),
because this supremum goes to infinity for each orbit that is not eventually
periodic (see \cite{WB}).

This difficulty may be overcome in two ways:

1) by considering open covers instead of partitions. This approach was
proposed by Brudno \cite{brud}. Since the sets in an open cover can have non
empty intersection, a step of the orbit of $x$ can be contained at time in
more than one open set of the cover. This implies that an orbit may have an
infinite family of possible symbolic codings, among which we choose the
``simplest one'';

2) by considering only a particular class of partitions (that are computable
in some sense that will be clarified later) and define the orbit complexity of
a point as the supremum of the orbit complexity over that class.

Brudno's open cover construction is not suitable for computational purposes
because the choice of the simplest coding in a big family is not practically feasible.

On the other hand, the computable partition approach is the mathematical
formalization of what we do in computer simulations. We consider a dynamical
system $(X,T)$ and we choose a partition $\beta$, which is always computable
when it is explicitly given (the formal definition will be given in next
section). We consider the symbolic orbit of a point $x\in X$ with respect to
$\beta$: it is a single string and we measure its information content by some
suitable universal coding algorithm.

\subsection{Computable partitions}

In this section we will give a rigorous definition of computable partition.
This notion is based on the idea of computable structure which relates the
abstract notion of metric space with computer simulations. Before giving the
formal definitions, few words are necessary to describe the intuitive idea
about it. Many models of the real words use the notion of real numbers or more
in general the notion of complete metric spaces. Even if you consider a very
simple complete metric space, as, for example, the interval $\left[
0,1\right]  $ it contains a continuum of elements. This fact implies that most
of these elements (numbers) cannot be described by any finite alphabet.
Nevertheless, in general, the mathematics of complete metric spaces is simpler
than the ''discrete mathematics'' in making models and the relative theorems.
On the other hand the discrete mathematics allows to make computer
simulations. A first connection between the two worlds is given by the theory
of approximation. But this connection becomes more delicate when we need to
simulate more sophisticated objects of continuum mathematics. For example an
open cover or a measurable partition of $\left[  0,1\right]  $ is very hard to
be simulated by computer; nevertheless, these notions are crucial in the
definition of many quantities as e. g. the K-S entropy of a dynamical system
or the Brudno complexity of an orbit. For this reason, we have introduced the
notion of ''computable structure'' which is a new way to relate the world of
continuous models with the world of computer simulations.

To simplify notations in the following, we will denote by $\Sigma$ the set of
finite binary strings in a finite alphabet. $\Sigma$ is the mathematical
abstraction of the world of the ''computer'', or more in general is the
mathematical abstraction of the ''things'' which can be expressed by any
language. We suppose that the real objects which we want to talk about are
modeled by the elements of a metric space $(X,d).$ We want to interpret the
objects of $\Sigma$ as points of $X.$ Thus, the first definition we need is
the following one:

\begin{definition}
An interpretation function on $(X,d)$ is a function $\Psi:\Sigma
_{0}\rightarrow X$ such that $\Psi(\Sigma)$ is dense in $X$. ($\Sigma
_{0}\subset\Sigma\;$is supposed to be a recursive set, namely a set for which
there is an algorithm which ''decides'' wether a given element belongs to set
or not).
\end{definition}

For example, if $X=\left[  0,1\right]  ,$ a possible choice of $\Psi$ is the
following one: if $s=s_{1}\dots s_{n}\in\Sigma=\Sigma_{0}=\Sigma(0,1),$ then
\begin{equation}
\Psi(s)=\sum_{i=1}^{n}s_{i}2^{-i}. \label{bubu}%
\end{equation}

A point $x\in X$ is said to be \emph{ideal} if it is the image of some string
$x=\Psi(s),\;s\in\Sigma$. Clearly almost every point is not ideal; however, it
is possible to perform computations on these points since they can be
approximated with arbitrary precision by ideal points, provided that the
interpretation is consistent with the distance $d.$ This consideration leads
us to the notion of ''computable interpretation'': an interpretation is
''computable'' if the distance between ideal points is computable with
arbitrary precision. The precise definition is the following:

\begin{definition}
A computable interpretation function on $(X,d)$ is a function $\Psi:\Sigma
_{0}\rightarrow X$ such that $\Psi(\Sigma_{0})$ is dense in $X$ and there
exists a total recursive function $D:\Sigma_{0}\times\Sigma_{0}\times
\rightarrow$ such that $\forall s_{1},s_{2}\in\Sigma,n\in\mathbf{N}$
\[
|d(\Psi(s_{1}),\Psi(s_{2}))-D(s_{1},s_{2},n)|\leq\frac{1}{2^{n}}.
\]
\end{definition}

If we take again $X=\left[  0,1\right]  $, we may have a different
interpretation considering a string in $\Sigma$ as an ASCII string. In this
case $\Sigma_{0}$ is the set of ASCII strings which denote a rational number
in $\left[  0,1\right]  .$ In this way we obtain a different computable
interpretation $\Psi:\Sigma_{0}\rightarrow X$ which describes the same metric
space. For most practical purposes, these two computable interpretations are
essentially equivalent: they represent the same \textit{computable structure.
}A computable structure on a separable metric space $(X,d)$ is a class of
\textit{equivalent }computable interpretations;.two interpretations are said
to be equivalent if the distance of ideal points is computable up to arbitrary precision.

\begin{definition}
\label{def7} Let $\Psi_{1}:\Sigma_{1}\rightarrow X$ and $\Psi_{2}:\Sigma
_{2}\rightarrow X$ be two computable interpretations in $(X,d)$; we say that
$\Psi_{1}$ and $\Psi_{2}$ are equivalent if there exists a total recursive
function $D^{\ast}:\Sigma_{1}\times\Sigma_{2}\times\mathbf{N}\rightarrow
\mathbf{R}$, such that $\forall s_{1},s_{2}\in\Sigma_{1}\times\Sigma
_{2},\;n{\in}\mathbf{N}$:
\[
|d(\Psi_{1}(s_{1}),\Psi_{2}(s_{2}))-D^{\ast}(s_{1},s_{2},n)|\leq\frac{1}%
{2^{n}}.
\]
\end{definition}

\begin{proposition}
The relation defined by definition \ref{def7} is an equivalence relation.
\end{proposition}

For the proof of this proposition see \cite{gal2}.

\begin{definition}
A computable structure $\mathcal{J}$ on $X$ is an equivalence class of
computable interpretations in $X$.
\end{definition}

Many concrete metric spaces used in analysis or in geometry have a natural
choice of a computable structure. The use of computable structures allows to
consider algorithms acting over metric spaces and to define constructive
functions between metric spaces (see \cite{gal2} or \cite{gal3}).

For example if $X=\mathbf{R}$ we can consider the interpretation $\Psi
:\Sigma\rightarrow\mathbf{R}$ defined in the following way: if $s=s_{1}\dots
s_{n}\in\Sigma$ then
\begin{equation}
\Psi(s)=\sum_{1\leq i\leq n}s_{i}2^{[n/2]-i}. \label{I}%
\end{equation}

This is an interpretation of a string as a binary expansion of a number.
$\Psi$ is a computable interpretation, the \emph{standard} computable
structure on $\mathbf{R.}$

\begin{definition}
Let $X$ be a space with a computable structure $\mathcal{I}$ and $\Psi
\in\mathcal{I}$. A partition $\beta=\{B_{i}\}$ is said to be computable if
there exists a recursive $Z:\Sigma\rightarrow$ such that for each $s$,
$Z(s)=j\Longleftrightarrow\Psi(s)\in B_{j}\ .$
\end{definition}

For example let us consider the following partition $\beta_{s_{0}%
,n}=\{B_{s_{0}}^{n},X-B_{s_{0}}^{n}\}$ with
\[
B_{s_{0}}^{n}=\overline{\{\Psi(s),D(s_{0},s,n+2)\leq\frac{1}{2^{n}}\}}%
-\{\Psi(s),D(s_{0},s,n+2)>\frac{1}{2^{n}}\}.
\]
\indent We remark that if $n$ goes to infinity, the diameter of the set $B_{s_{0}}%
^{n}$ goes to 0. Moreover since $X$ is compact, for each $n$ there is a finite
set of strings $S_{n}=\{s_{0},\dots,s_{k}\}$ such that $X=\cup_{s\in S_{n}%
}B_{s}^{n}$. Since the join of a finite family of computable partitions is
computable, there is a family of partitions $\alpha_{n}=\vee_{s\in S_{n}}%
\beta_{s,n}$ such that $\alpha_{n}$ is for each $n$ a finite computable
partition and $lim_{n\rightarrow\infty}diam(\alpha_{n})=0$. This will be used
in the proof of the next theorem.

\subsection{Computable complexity of an orbit}

Using the notion of computable partition, it is possible to define the
complexity of a single orbit just using notion given by information theory.
Since we have to different notions of complexity of a srting (namely
$K_{AIC}(s)$ and $K_{Z}(s)$), we get two different notions of complexity of an orbit:

\begin{definition}
If $(X,T)$ is a dynamical system over a metric space $(X,d)$ with a computable
structure we define the computable complexity of the orbit of $x$ as:
\[
K_{AIC}(x)=\sup\{K_{AIC}(x,\beta)\ |\ \beta\ \mbox{computable partition}\}
\]%
\[
K_{Z}(x)=\sup\{K_{Z}(x,\beta)\ |\ \beta\ \mbox{computable partition}\}.
\]
\end{definition}

\begin{theorem}
If $(X,\mu,T)$ is a dynamical system on a compact space and $\mu$ is ergodic,
then for $\mu$-almost each $x,$
\[
K_{Z}(x)=K_{AIC}(x)=h_{\mu}(T)
\]
\end{theorem}

\emph{Proof.} By what it is said above we remark that for each $\epsilon$
there is a computable partition with diameter less than $\epsilon$. Since
computable partitions are a countable set, by lemma \ref{15} we prove the statement.
$\Box$

The above theorem states that $K_{AIC}(x)$ and $K_{Z}(x)$ are the right
quantity to be considered; in fact in dynamical systems which are stationary
and ergodic they coincide with the K-S entropy for $a.e.$ $x$. However the
basic point is that $K_{Z}(x),$in principle, can be computed by a machine with
an arbitrary accuracy.

We call the function $K_{Z}(x)$ \emph{computable complexity} of an orbit. For
stationary and ergodic dynamical systems, it is equivalent to the K-S entropy;
nevertheless, it has a broader applicability and it presents the following features:

\begin{itemize}
\item  it is the limit of computable functions and it can be measured in real experiments;

\item it is independent of the invariant measure;

\item  it makes sense even if the space $X$ is not compact;

\item  it makes sense even if the dynamics is not stationary.
\end{itemize}

\section{Numerical Experiments\label{num}}

Weakly chaotic dynamical systems give symbolic orbits with null entropy. For
these systems the behavior of the quantity of information that is contained in
$n$ steps of the symbolic orbit is less than linear. However there is a big
difference between a periodic dynamical system and a sporadic one (for example
the Manneville map, Section 4.1). In fact, the latter can have positive
topological entropy and sensitive dependence on initial conditions.

Thus it is important to have a numerical way to detect weak chaos and to
classify the various kind of chaotic behavior.

We have implemented a particular compression algorithm which we called CASToRe
(Compression Algorithm, Sensitive To Regularity). Here we present some
experiments. First we used CASToRe in the study of the Manneville map. The
information behavior of the Manneville map is known by the works of
\cite{Gaspard}, \cite{gal3}, (where the $AIC$ has been used), \cite{bonanno} .
We will see that our me\-thods, implemented by a computer simulation, give
results which coincide with the the theoretical predictions of the mentioned
papers. In the second example, CASToRe is used to have a measure of the kind of
chaotic behavior of the logistic map at the chaos threshold. Previous
numerical works from the physics literature proved that the logistic map has
power law data sen\-si\-ti\-vi\-ty to initial conditions, which implies logarithmic
growth of the quantity of information. Our numerical results surely suggest
that in such a dynamical system the behavior of the quantity of information is
below any power law, confirming the previous results.

\subsection{The Manneville map}

\label{manne} 
\begin{figure}[th]
\centerline{
\psfig{figure=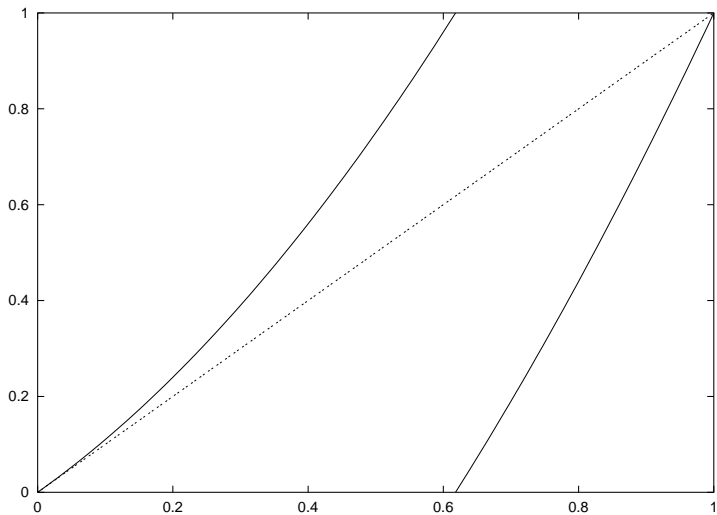,width=10cm,angle=270}}\caption{\textit{The Manneville
map $f$ for $z=2$}}
\label{grafmann}
\end{figure}

The \textit{Manneville map} was introduced by Manneville in
\cite{Manneville} as an example of a discrete dissipative dynamical system
with \textit{intermittency}, an alternation between long regular phases,
called \textit{laminar}, and short irregular phases, called \textit{turbulent}%
. This behavior has been observed in fluid dynamics experiments and in
chemical reactions. Manneville introduced his map, defined on the interval
$I=[0,1]$ by
\begin{equation}
f(x)=x+x^{z}(\hbox{mod }1)\hskip0.5cmz>1, \label{mannmap}%
\end{equation}
to have a simple model displaying this complicated behavior (see Figure
\ref{grafmann}). His work has attracted much attention, and the dynamics of
the Manneville map has been found in many other systems. We can find
applications of the Manneville map in dynamical approaches to DNA sequences
(\cite{Grigolini},\cite{Grigolini96}) and ion channels (\cite{Toth}), and in
non-extensive thermodynamic problems (\cite{Grigo3}).

Gaspard and Wang (\cite{Gaspard}) described the behavior of the algorithmic
complexity of the Manneville map, and called such a behavior
\textit{sporadicity}. A mathematical study of this kind of behavior has been
done in a more general contest (\cite{bonanno}) and using different methods
(\cite{gal3}).

Let's consider the following partition of the interval $[0,1]$. Let $x_{0}$ be
such that $f(x_{0})=1$ and $x_{0}\not =1$ (see Figure \ref{grafmann}), then we
call $A_{0}=(x_{0},1]$ and $A_{1}=[0,x_{0}]$. We can now recursively define
points $x_{k}$ such that $f(x_{k})=x_{k-1}$ and $x_{k}<x_{k-1}$. Then we
consider the partition $\alpha=(A_{k})$ where $A_{k}=(x_{k},x_{k-1})$, for
$k\in$. We denote by $I_{AIC}(x,\alpha,n)$ the Algorithmic Information Content
of a $n$-long orbit of the Manneville map with initial condition $x$, using
the partition $\alpha$. We have that the mean value of $I_{AIC}(x,\alpha,n)$,
with respect to the Lebesgue measure $l$, on the initial conditions of the
orbit is $E_{l}[I_{AIC}(x,\alpha,n)]\sim n^{p}$, with $p=\frac{1}{z-1}$ for
$z\geq2$, and $E_{l}[I_{AIC}(x,\alpha,n)]\sim n$ for $z<2$.

Experiments which have been performed by using the compression algorithm CASToRe (\cite{ZIPPO})
confirm the theoretical results and prove that the method related to the
computable complexity are experimantally reliable. We considered a set of one
hundred initial points, generated $10^{7}$-long orbits, and applied the
algorithm to the associated symbolic strings $s$. If we considered the
compression algorithm $Z=$CASToRe, we have that $I_{Z}(s)$ is a good
approximation of $I_{AIC}(x,\alpha,n)$.

In Table \ref{tabellamann} we show the results. The first column is the value
of the parameter $z$. The last column gives the results of the theory for the
exponent of the asymptotic behavior of $I_{AIC}(x,\alpha,n)$. The second and
third column show the experimental results. Given the functions $I_{Z}(s)$,
with $Z=$CASToRe, we made a mean of the obtained values for the exponents $p$
using the Lebesgue measure (second column) and the invariant density given by
the system (third column) (\cite{Gaspard}).

\begin{table}[th]%
\begin{tabular}
[c]{|c|c|c|c|}\hline
z & uniform distribution & invariant distribution & theoretical value\\\hline
2.8 & 0.573 & 0.551 & 0.555\\\hline
3 & 0.533 & 0.509 & 0.5\\\hline
3.5 & 0.468 & 0.449 & 0.4\\\hline
4 & 0.409 & 0.381 & 0.333\\\hline
\end{tabular}
\caption{Theoretical and experimental results for the Manneville map}%
\label{tabellamann}%
\end{table}

\begin{figure}[ptb]%
\begin{tabular}
[c]{lr}%
{\raggedright{\psfig{figure=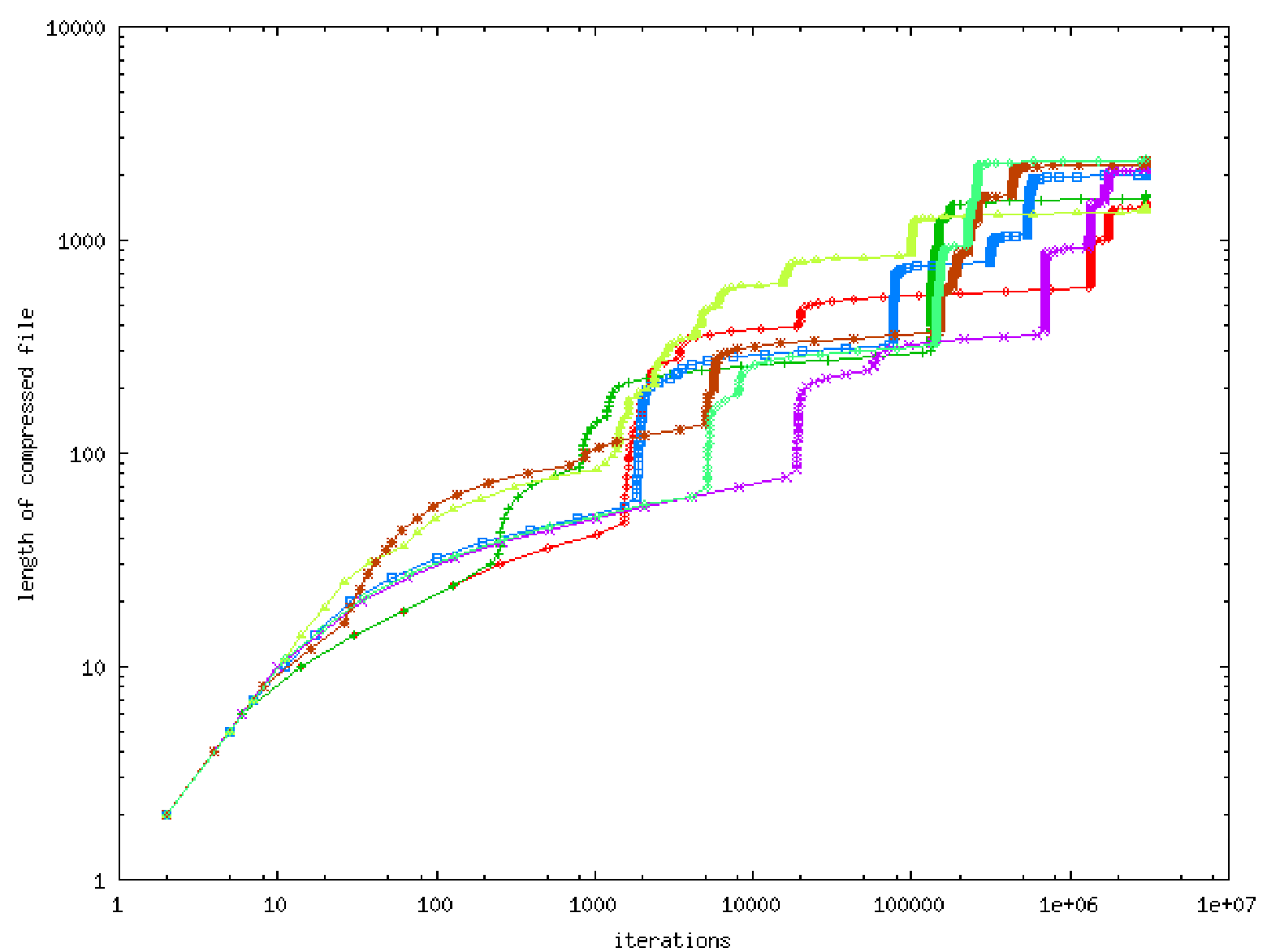,width=5.5cm,angle=0}}} & {\raggedleft{
\psfig{figure=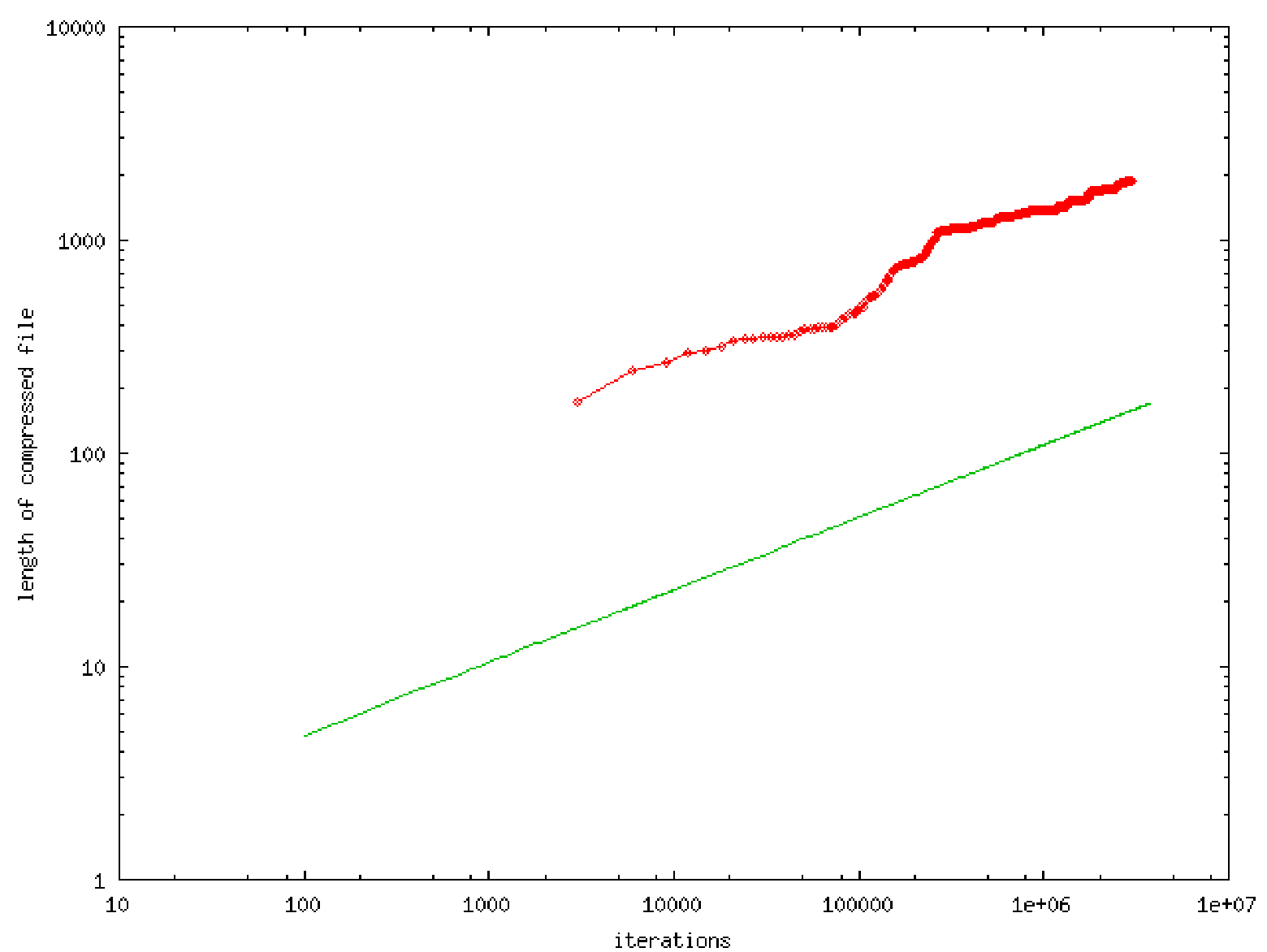,width=5.5cm,angle=0}}}%
\end{tabular}
\caption{\textit{The experimental results for seven different orbits generated
by the Manneville map with $z=4$}}%
\label{figmann}%
\end{figure}

These experiments seem also to show, as it is indicated by the theory
(\cite{bonanno}), that the Algorithmic Information Content $I_{AIC}%
(x,\alpha,n)$ of strings generated by the Manneville map is such that
$I_{AIC}(x,\alpha,n) \sim E_{l}[I_{AIC}(x,\alpha,n)]$ for almost any initial
condition with respect to the Lebesgue measure $l$. In Figure \ref{figmann},
we show the experimental results for the Manneville map with $z=4$. On the
left there are plotted the functions $I_{Z}(s)$, with $Z=$CASToRe, for seven
different initial conditions, and on the right there is the mean of the
functions and a right line with slope $0.33$, showing the asymptotic behaviour
$n^{p}$. Notice that the functions are plotted with logarithmic coordinates,
then $p$ is the slope of the lines. 

\subsection{The logistic map\label{logistic}}

We studied the algorithmic complexity also for the logistic map defined by
\begin{equation}
f(x)=\lambda x(1-x)\ ,\quad x\in\lbrack0,1]\ ,\quad1\leq\lambda\leq4.
\label{logmap}%
\end{equation}

The logistic map has been used to simulate the behavior of biological species
not in competition with other species. Later the logistic map has also been
presented as the first example of a relatively simple map with an extremely
rich dynamics. If we let the parameter $\lambda$ vary from $1$ to $4$, we find
a sequence of bifurcations of different kinds. For values of $\lambda
<\lambda_{\infty}=3.56994567187\dots$, the dynamics is periodic and there is a
sequence of period doubling bifurcations which leads to the chaos threshold
for $\lambda=\lambda_{\infty}$. The dynamics of the logistic map at the chaos
threshold has attracted much attention and many are the applications of
theories to the logistic map at this particular value of the parameter
$\lambda$. In particular, numerical experiments suggested that at the chaos
threshold the logistic map has null K-S entropy, implying that the sensitivity
on initial conditions is less than exponential, and there is a power-law
sensitivity to initial conditions. These facts have justified the application
of generalized entropies to the map (\cite{tsallis1},\cite{tsallis}).
Moreover, from the relations between initial conditions sensitivity and
information content (\cite{gal3}), we expect to find that the Algorithmic
Information Content of the logistic map at the chaos threshold is such that
$I_{AIC}(s)\sim\log(n)$ for a $n$-digit long symbolic string $s$ generated by
one orbit of the map. We next show how we have experimentally found this result.

\begin{figure}[ptb]
\centerline{\psfig{figure=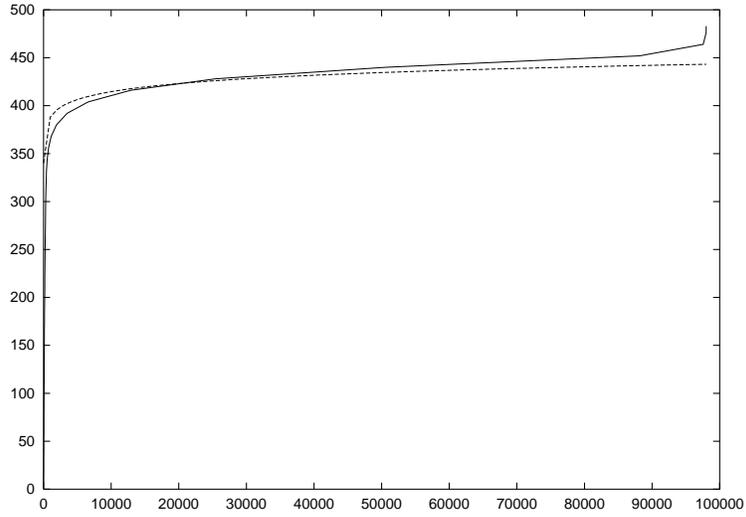,width=10cm,angle=270}}\caption{\textit{The solid line
is the information function for a string $\sigma$ of period 100 compared with
$\Pi(N)$ plus a constant term}}
\label{figperlog}
\end{figure}

From now on $Z$ indicates the compression algorithm CASToRe. It is known that
for periodic maps the behavior of the Algorithmic Information Content should
be of order $O(\log n)$ for a $n$-long string, and it has been proved that the
compression algorithm CASToRe gives for periodic strings $I_{Z}(s)\sim
\Pi(n)=\log n\log(\log n)$, where we recall that $I_{Z}(s)$ is the binary
length of the compressed string (\cite{menconi}). In Figure \ref{figperlog},
we show the approximation of $I_{Z}(s)$ with $\Pi(n)$ for a $10^{5}$-long
periodic string of period 100. 

We have thus used the sequence $\lambda_{n}$ of parameters values where the
period doubling bifurcations occur and used a ``continuity'' argument to
obtain the behavior of the information function at the chaos threshold.
Another sequence of parameters values $\mu_{k}$ approximating the critic value
$\lambda_{\infty}$ from above has been used to confirm the results.

\begin{figure}[ptb]
\centerline{\psfig{figure=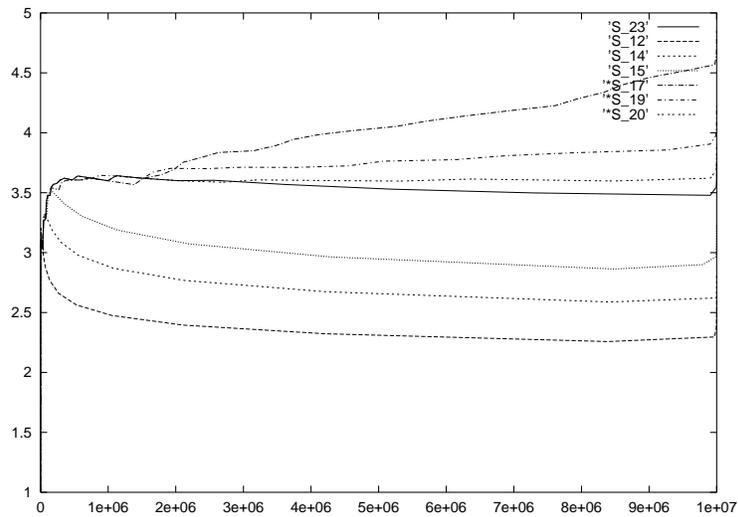,width=10cm,angle=270}}\caption{\textit{The solid line is
the limit function $S_{\infty}(n)$ and dashed lines are some approximating
functions from above (with a star) and below.}}%
\label{figapproxlog}%
\end{figure}

In Figure \ref{figapproxlog}, we plotted the functions $S(n)=\frac{I_{Z}%
(s)}{\Pi(n)}$ for some values of the two sequences. The starred functions
refer to the sequence $\mu_{k}$ and the others to the sequence $\lambda_{n}$.
The solid line show the limit function $S_{\infty}(n)$. If we now consider the
limit for $n\rightarrow+\infty$, we conjecture that $S_{\infty}(n)$ converges
to a constant $S_{\infty}$, whose value is more or less $3.5$. Then we can
conclude that, at the chaos threshold, the Algorithmic Information Content of
the logistic map is $I_{AIC}(s)\leq I_{Z}(s)\sim S_{\infty}\Pi(n)$. In
particular we notice that we obtained an Algorithmic Information Content whose
order is smaller than any power law, and we called this behavior \textit{mild
chaos} (\cite{menconi}). 

\section{DNA sequences\label{DNAm}}

\label{DNA} We look at genomes as to finite symbolic sequences where the
alphabet is the one of nucleotides $\{\mbox{A, C, G, T}\}$ and compute the
complexity $K_{Z}$ (where $Z$ is the algorithm CASToRe) of some sequences or
part of them.

DNA sequences, in fact, can be roughly divided in different functional
regions. First, let us analyze the structure of a Prokaryotic genome: a gene
is not directly translated into protein, but is expressed via the production
of a messenger RNA. It includes a sequence of nucleotides that corresponds
exactly with the sequence of amino acids in the protein (this is the so called
\textit{colinearity} of prokaryotic genes). These parts of the genome are the
\textbf{coding regions}. The regions of prokaryotic genome that are not coding
are the \textbf{non coding regions}: upstream and downstream regions, if they
are proceeding or following the gene.

On the other hand, Eukaryotic DNA sequences have several non coding regions: a
gene includes additional sequences that lie within the coding region,
interrupting the sequence that represent the protein (this is why these are
\textit{interrupted genes}). The sequences of DNA comprising an interrupted
gene are divided into two categories:

\begin{itemize}
\item  the \textbf{exons} are the sequences represented in the mature RNA and
they are the real coding region of the gene, that starts and ends with exons;

\item  the \textbf{introns} are the intervening sequences which are removed
when the first transcription occurs.
\end{itemize}

So, the non coding regions in Eukaryotic genomes are intron sequences and
up/downstream sequences. The last two regions are usually called
\textbf{intergenic regions}. In Bacteria and Viruses genomes, coding regions
have more extent than in Eukaryotic genomes, where non coding regions prevail.

There is a long-standing interest in understanding the correlation structure
between bases in DNA sequences. Statistical heterogeneity has been
investigated separately in coding and non coding regions: long-range
correlations were proved to exist in intron and even more in intergenic
regions, while exons were almost indistinguishable from random sequences
(\cite{Li vari},\cite{AcquistiBuiatti}). \begin{figure}[ptb]
\begin{tabular}
[c]{lr}
{\raggedright{\psfig{figure=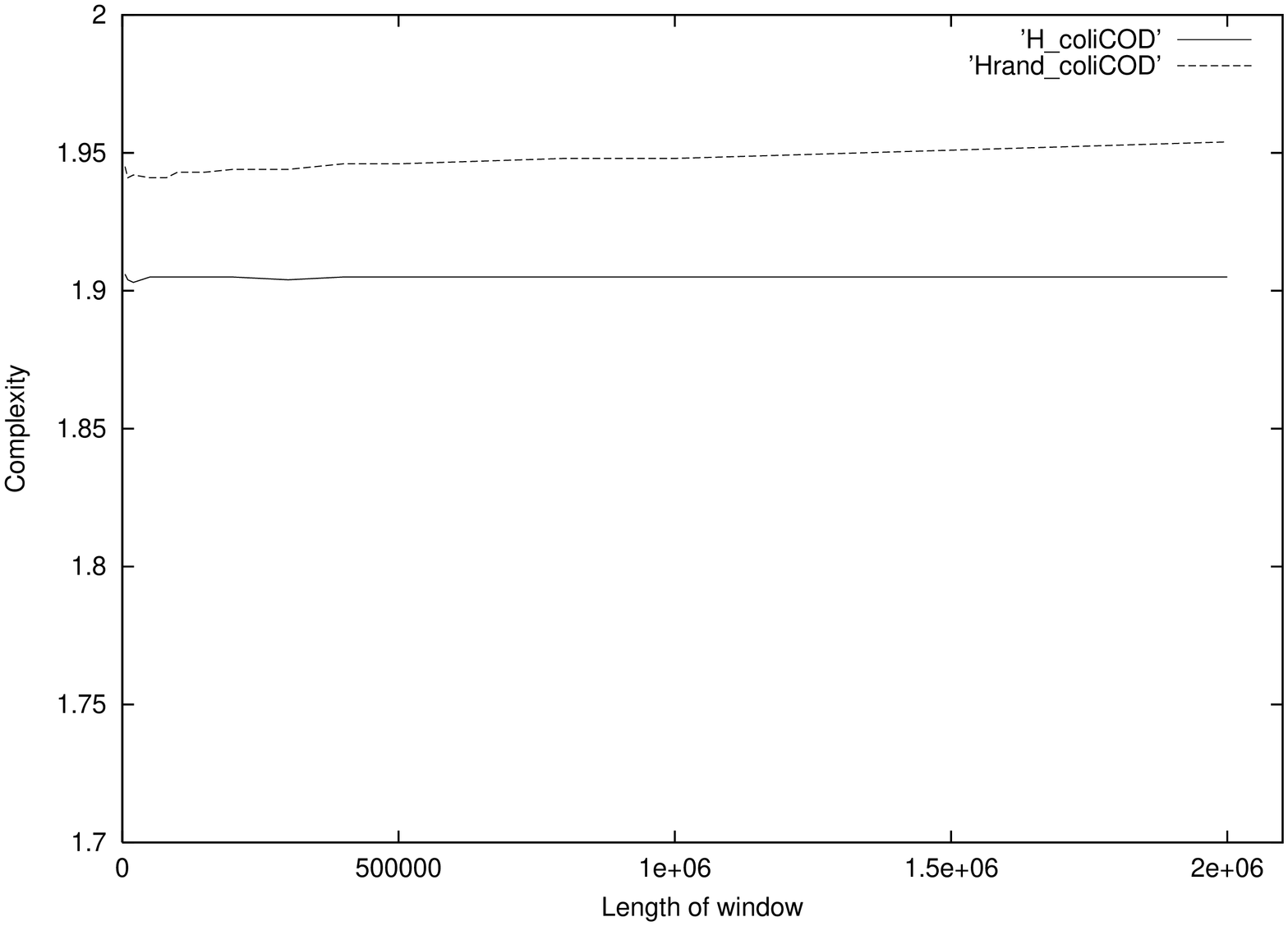,width=6cm,angle=270}}} &
{\raggedleft{\psfig{figure=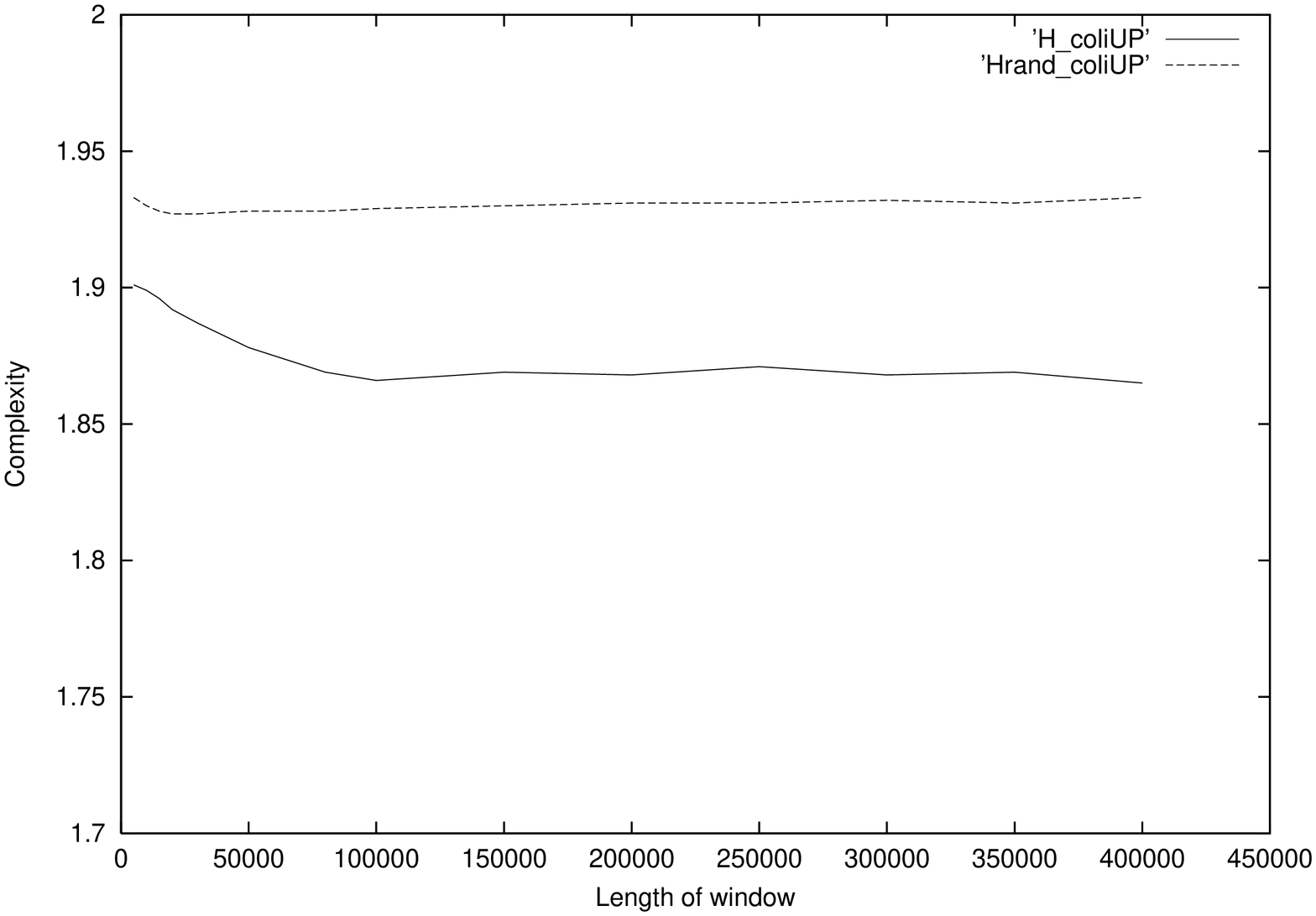,width=6cm,angle=270}}}
\end{tabular}
\caption{\textit{On the left, the complexity of coding sequences of Bacterium
Escherichia Coli as a function of the length $L$ of the windows are compared
to a statistically equivalent random sequence. On the right, the same
comparison with respect to the intergenic (non coding) sequence of the same
genome.}}
\label{coli+rand}
\end{figure}\begin{figure}[ptbptb]
\centerline{
\psfig{figure=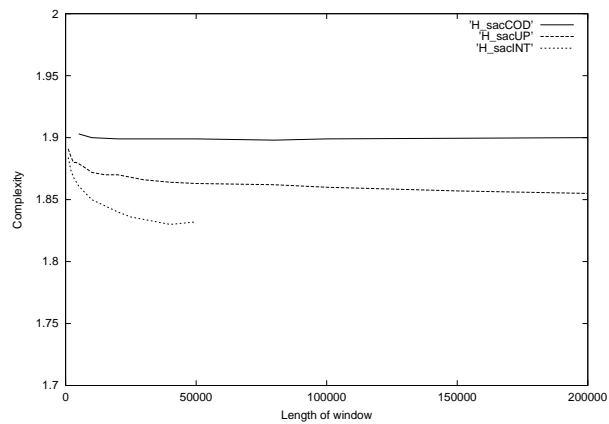,width=8cm,angle=270}}
\caption{\textit{The analysis of the three functional regions of the genome of
Eukaryote Saccharomyces Cerevisiae: the solid line is relative to the coding
regions, the dashed line to the intergenic regions, the dotted line to the
intron region. The complexities, as functions of length $L$ of the windows,
are clearly different from each other. The regions have not similar length so
the picture has been drawn according to the shortest region.}}%
\label{cfrSAC}%
\end{figure}Our approach can be applied to look for the
\textbf{non-extensivity of the Information content corresponding to the
different regions of the sequences}.

We have used a modified version of the algorithm CASToRe, that exploits a
window segmentation (see Appendix); let $L$ be the length of a window. We
measure the mean complexity of substrings with length $L$ belonging to the
sequence that is under analysis. Then, we obtain the Information and the
complexity of the sequence as functions of the length $L$ of the windows.

In analogy with physical language, we call a function $f$ extensive if the sum
of the evaluations of $f$ on each part of a partition of the domain equals the
evaluation of $f$ on the whole domain:
\[
\sum_{i\in I}f(A_{i})=f\left(  \bigcup_{i\in I}A_{i}\right)  .
\]

In case of a non-extensive function $\tilde{f}$, we have that the average on
the different parts underestimates the evaluation on the whole domain:
\[
\sum_{i\in I}\tilde{f}(A_{i})\leq\tilde{f}\left(  \bigcup_{i\in I}%
A_{i}\right)  .
\]

Now let us consider the Information Content: each set $A_{i}$ is a window of
fixed length $L$ in the genome, so it can be considered as $A_{i}(L)$. Then
the related complexity is $\frac{I_{Z}(A_{i}(L))}{L}$.

If a sequence is chaotic or random, then its Information content is extensive,
because it has no memory (neither short-range or long range correlations). But
if a sequence shows correlations, then the more long-ranged the correlations
are the more non-extensive is the related Information content. We have that:

\begin{itemize}
\item  if the Information content is extensive, the complexity is constant as
a linear function of the length $L$;

\item  if the Information content is non-extensive, the complexity is a
decreasing, less than linear function of the length $L$.
\end{itemize}

From the experimental point of view, we expect our results to show that in
coding regions the Information content is extensive, while in non coding
regions the extensivity is lost within a certain range $[0,L^{\ast}]$ of
window length (the number $L^{\ast}$ depends on the genome). This is also
supported by the statistical results exposed above.

In coding sequences, we found that the complexity is almost constant: the
compression ratio does not change sensitively with $L$. In non coding
sequences, the complexity decreases until some appropriate $L^{\ast}$ is
reached and later the compression ratio is almost constant.

This is an information-theoretical proof (alternative to the statistical
technique) that coding sequences are more chaotic than non coding ones. Figure
\ref{coli+rand} shows the complexity of coding (on the left) and non coding
(on the right) regions of the genome of Bacterium Escherichia Coli as a
function of the length $L$ of the windows and compared with statistically
equivalent random sequences. Clearly the compression ratio decreases more in
the non coding regions than in the coding ones, where the random sequences
show almost constant complexity.

Figure \ref{cfrSAC} shows the analysis of the three functional regions of the
genome of Saccharomyces Cerevisiae which is an eukaryote: the complexities, as
functions of length $L$ of the windows, are clearly different from each other.
We remark that the lower is the compression ratio, the higher is the
aggregation of the words in the genome. This is due to the fact that the
algorithm CASToRe recognizes patterns already read, so it is not necessary to
introduce new words, but coupling old words to give rise to a new longer one
is sufficient to encode the string.

\section*{Appendix: CASToRe}

\label{appendice}
The Lempel-Ziv $LZ78$ coding scheme is traditionally used to codify a
string according to a so called incremental parsing procedure
\cite{LZ78}. The algorithm divides the sequence in words that are all
different from each other and whose union is the dictionary. A new
word is the longest word in the charstream that is representable as a
word already belonging to the dictionary together with one symbol
(that is the ending suffix). 

We remark that the algorithm $LZ78$
encodes a constant $n$ digits long sequence $'111\dots'$ to a string
with length about $const\ +\ n^{\frac 1 2}$ bits, while the
theoretical Information $I_{AIC}$ is about $const\ +\ \log n$. So, we
can not expect that $LZ78$ is able to distinguish a sequence whose
Information grows like $n^{\alpha}$ ($\alpha < \frac 1 2$) from a
constant or periodic one. 

This is the main motivation which lead us to
create the new algorithm CASToRe. It has been proved in \cite{menconi}
that the Information of a constant sequence, originally with length
$n$, is $4+2\log (n+1)[\log (\log (n+1))-1]$, if CASToRe is used. As
it has been showed in section 4.1, the new algorithm is also a
sensitive device to weak
dynamics.

CASToRe is an encoding algorithm based on an adaptive dictionary. Roughly
speaking, this means that it translates an input stream of symbols (the file
we want to compress) into an output stream of numbers, and that it is possible
to reconstruct the input stream knowing the corrispondence between output and
input symbols. This unique corrispondence between sequences of symbols (words)
and numbers is called \emph{the dictionary}.

``Adaptive'' means that the dictionary depends on the file under compression,
in this case the dictionary is created while the symbols are translated.

At the beginning of encoding procedure, the dictionary is empty. In order to
explain the principle of encoding, let's consider a point within the encoding
process, when the dictionary already contains some words.

We start analyzing the stream, looking for the longest word W in the
dictionary matching the stream. Then we look for the longest word Y in the
dictionary where W + Y matches the stream. Suppose that we are compressing an
english text, and the stream contains ``basketball ...'', we may have the
words ``basket'' (number 119) and ``ball'' (number 12) already in the
dictionary, and they would of course match the stream.

The output from this algorithm is a sequence of word-word pairs (W, Y), or
better their numbers in the dictionary, in our case (119, 12). The resulting
word ``basketball'' is then added to the dictionary, so each time a pair is
output to the codestream, the string from the dictionary corresponding to W is
extended with the word Y and the resulting string is added to the dictionary.

A special case occurs if the dictionary doesn't contain even the starting
one-character string (for example, this always happens in the first encoding
step). In this case we output a special code word which represents the null
symbol, followed by this character and add this character to the dictionary.

Below there is an example of encoding, where the pair (4, 3), is composed from
the fourth word $'$AC$'$ and the third word $'$G$'$.

ACGACACGGAC

\ 

word 1: $(0, A)$ --$>$ A

word 2: $(0, C)$ --$>$ C

word 3: $(0, G)$ --$>$ G

word 4: $(1, 2)$ --$>$ AC

word 5: $(4, 3)$ --$>$ ACG

word 6: $(3, 4)$ --$>$ GAC

This algorithm can be used in the study of correlations.

A modified version of the program can be used for study of correlations in the
stream: it partition the stream into fixed size segments and proceed to encode
them separately. The algorithm takes advantage of replicated parts in the
stream. Limiting the encoding to each window separately, could results in
longer total encoding. The difference between the length of the whole stream
encoded and the sum of the encoding of each window depends on number of
correlation between symbols at distance greater than the size of the window.
Using different window sizes it is possible to construct a "spectrum" of correlation.

This has been applied for example on DNA sequences for the study of mid and
long-term correlations (see section \ref{DNAm}).
\newpage
\textbf{Implementation:}

The main problem implementing this algorithm is building a structure which
allows to efficently search words in the dictionary. To this purpose, the
dictionary is stored in a treelike structure where each node is a word X = (W,
Y), its parent node being W, and storing a link to the string representing Y.
Using this method, in order to find the longest word in the dictionary maching
exactly a string, you need only to follow a branching path from the root of
the dictionary tree.

\end{document}